\documentclass{article}

\usepackage{amssymb,amsmath,theorem,euscript}

\newcounter{sec}

\newcounter{punct}[sec]

\def\punct{\refstepcounter{punct}{\arabic{sec}.\arabic{punct}.  }}

\def\COUNTERS{\addtocounter{sec}{1}
              \setcounter{punct}{0}
          \setcounter{equation}{0}
          \setcounter{theorem}{0}
                  }

\newtheorem{theorem}{Theorem}[sec]
\newtheorem{proposition}[theorem]{Proposition}

\newtheorem{lemma}[theorem]{Lemma}

 \def\ov{\overline}
\def\wt{\widetilde}

\begin{document}

\def\OO{\mathrm{O}}
\def\GLO{\mathrm{GLO}}
\def\Coll{\mathrm{Coll}}
\def\kappa{\varkappa}

\def\R{\mathbb{R}}
\def\C{\mathbb{C}}

\def\la{\langle}
\def\ra{\rangle}

 \def\cA{\mathcal A}
\def\cB{\mathcal B}
\def\cC{\mathcal C}
\def\cD{\mathcal D}
\def\cE{\mathcal E}
\def\cF{\mathcal F}
\def\cG{\mathcal G}
\def\cH{\mathcal H}
\def\cJ{\mathcal J}
\def\cI{\mathcal I}
\def\cK{\mathcal K}
 \def\cL{\mathcal L}
\def\cM{\mathcal M}
\def\cN{\mathcal N}
 \def\cO{\mathcal O}
\def\cP{\mathcal P}
\def\cQ{\mathcal Q}
\def\cR{\mathcal R}
\def\cS{\mathcal S}
\def\cT{\mathcal T}
\def\cU{\mathcal U}
\def\cV{\mathcal V}
 \def\cW{\mathcal W}
\def\cX{\mathcal X}
 \def\cY{\mathcal Y}
 \def\cZ{\mathcal Z}
\def\0{{\ov 0}}
 \def\1{{\ov 1}}
 \def\frA{\mathfrak A}
 \def\frB{\mathfrak B}
\def\frC{\mathfrak C}
\def\frD{\mathfrak D}
\def\frE{\mathfrak E}
\def\frF{\mathfrak F}
\def\frG{\mathfrak G}
\def\frH{\mathfrak H}
\def\frI{\mathfrak I}
 \def\frJ{\mathfrak J}
 \def\frK{\mathfrak K}
 \def\frL{\mathfrak L}
\def\frM{\mathfrak M}
 \def\frN{\mathfrak N} \def\frO{\mathfrak O} \def\frP{\mathfrak P} \def\frQ{\mathfrak Q} \def\frR{\mathfrak R}
 \def\frS{\mathfrak S} \def\frT{\mathfrak T} \def\frU{\mathfrak U} \def\frV{\mathfrak V} \def\frW{\mathfrak W}
 \def\frX{\mathfrak X} \def\frY{\mathfrak Y} \def\frZ{\mathfrak Z} \def\fra{\mathfrak a} \def\frb{\mathfrak b}
 \def\frc{\mathfrak c} \def\frd{\mathfrak d} \def\fre{\mathfrak e} \def\frf{\mathfrak f} \def\frg{\mathfrak g}
 \def\frh{\mathfrak h} \def\fri{\mathfrak i} \def\frj{\mathfrak j} \def\frk{\mathfrak k} \def\frl{\mathfrak l}
 \def\frm{\mathfrak m} \def\frn{\mathfrak n} \def\fro{\mathfrak o} \def\frp{\mathfrak p} \def\frq{\mathfrak q}
 \def\frr{\mathfrak r} \def\frs{\mathfrak s} \def\frt{\mathfrak t} \def\fru{\mathfrak u} \def\frv{\mathfrak v}
 \def\frw{\mathfrak w} \def\frx{\mathfrak x} \def\fry{\mathfrak y} \def\frz{\mathfrak z} \def\frsp{\mathfrak{sp}}
 \def\bfa{\mathbf a} \def\bfb{\mathbf b} \def\bfc{\mathbf c} \def\bfd{\mathbf d} \def\bfe{\mathbf e} \def\bff{\mathbf f}
 \def\bfg{\mathbf g} \def\bfh{\mathbf h} \def\bfi{\mathbf i} \def\bfj{\mathbf j} \def\bfk{\mathbf k} \def\bfl{\mathbf l}
 \def\bfm{\mathbf m} \def\bfn{\mathbf n} \def\bfo{\mathbf o} \def\bfp{\mathbf p} \def\bfq{\mathbf q} \def\bfr{\mathbf r}
 \def\bfs{\mathbf s} \def\bft{\mathbf t} \def\bfu{\mathbf u} \def\bfv{\mathbf v} \def\bfw{\mathbf w} \def\bfx{\mathbf x}
 \def\bfy{\mathbf y} \def\bfz{\mathbf z} \def\bfA{\mathbf A} \def\bfB{\mathbf B} \def\bfC{\mathbf C} \def\bfD{\mathbf D}
 \def\bfE{\mathbf E} \def\bfF{\mathbf F} \def\bfG{\mathbf G} \def\bfH{\mathbf H} \def\bfI{\mathbf I} \def\bfJ{\mathbf J}
 \def\bfK{\mathbf K} \def\bfL{\mathbf L} \def\bfM{\mathbf M} \def\bfN{\mathbf N} \def\bfO{\mathbf O} \def\bfP{\mathbf P}
 \def\bfQ{\mathbf Q} \def\bfR{\mathbf R} \def\bfS{\mathbf S} \def\bfT{\mathbf T} \def\bfU{\mathbf U} \def\bfV{\mathbf V}
 \def\bfW{\mathbf W} \def\bfX{\mathbf X} \def\bfY{\mathbf Y} \def\bfZ{\mathbf Z} \def\bfw{\mathbf w}
 \def\R {{\mathbb R }} \def\C {{\mathbb C }} \def\Z{{\mathbb Z}} \def\H{{\mathbb H}} \def\K{{\mathbb K}}
 \def\N{{\mathbb N}} \def\Q{{\mathbb Q}} \def\A{{\mathbb A}} \def\T{\mathbb T} \def\P{\mathbb P} \def\G{\mathbb G}
 \def\bbA{\mathbb A} \def\bbB{\mathbb B} \def\bbD{\mathbb D} \def\bbE{\mathbb E} \def\bbF{\mathbb F} \def\bbG{\mathbb G}
 \def\bbI{\mathbb I} \def\bbJ{\mathbb J} \def\bbL{\mathbb L} \def\bbM{\mathbb M} \def\bbN{\mathbb N} \def\bbO{\mathbb O}
 \def\bbP{\mathbb P} \def\bbQ{\mathbb Q} \def\bbS{\mathbb S} \def\bbT{\mathbb T} \def\bbU{\mathbb U} \def\bbV{\mathbb V}
 \def\bbW{\mathbb W} \def\bbX{\mathbb X} \def\bbY{\mathbb Y} \def\kappa{\varkappa} \def\epsilon{\varepsilon}
 \def\phi{\varphi} \def\le{\leqslant} \def\ge{\geqslant}

\def\Gms{\mathrm{Gms}}
\def\Pol{\mathrm{Pol}}
\def\Mat{\mathrm{Mat}}
\def\vol{\mathrm{vol}}
\def\Lat{\mathrm{Lat}}

\def\U{\mathrm{U}}
\def\GL{\mathrm{GL}}
\def\SO{\mathrm{SO}}
\def\Sp{\mathrm{Sp}}
\def\OO{\mathrm{O}}
\def\Gr{\mathrm{Gr}}

\def\Symm{\mathrm{Symm}}
\def\ASymm{\mathrm{ASymm}}

\def\zigzag{\rightsquigarrow}

\def\sm{\smallskip}

\begin{center}
\Large\bf
Hua measures on the space of $p$-adic matrices and inverse limits of Grassmannians

\bigskip

\large\sc
Yury Neretin%
\footnote{Supported by grant FWF, P22122.}

\end{center}

{\small We construct $p$-adic counterparts of Hua measures, measures on
inverse limits of $p$-adic Grassmannians, and describe natural groups of symmetries 
of such measures.}

\section{Results of the paper}

\COUNTERS

{\bf\punct Real archetype.} For details, see \cite{Ner-hua}, \cite{Ner-gauss},
 Section II.2.10.  Denote by $\U(k)$ the group of unitary matrices of size $k$, by
 $d\chi^k$ the probability Haar measure on $\U(k)$.
By $\U(\infty)$ we denote the inductive limit
 $$
 \dots\longrightarrow \U(n)\longrightarrow \U(n+1)\longrightarrow\dots
 $$
 of groups $\U(n)$.
 
Represent an element $g\in\U(n+1)$ as a block matrix $g=\begin{pmatrix}a&b\\c&d \end{pmatrix}$
of size $n+1$. Consider
the {\it Livshits map} $\Upsilon:\U(n+1)\to\U(n)$ given by
\begin{equation}
\Upsilon \begin{pmatrix}a&b\\c&d \end{pmatrix}=a-b(1+d)^{-1}c
.
\label{eq:upsilon}
\end{equation}
This map commutes with left and right actions of $\U(n)$,
$$
\Upsilon(h_1 g h_2)=h_1\Upsilon(g) h_2, \qquad\text{where $h_1$, $h_2\in \U(n)$.}
$$
 Therefore, the pushforward
of the Haar measure $d\chi^{n+1}$ under the map $\Upsilon$  is the Haar measure
$d\chi^n$. Therefore there exists the inverse limit $\frU(\infty)$ of the chain
\begin{equation}
\dots\longleftarrow \U(n)\longleftarrow \U(n+1)\longleftarrow\dots
\label{eq:chain}
\end{equation}
equipped with a probability  measure $d\chi^\infty$. The space $\frU(\infty)$ is not a group,
 but the
unitary group $\U(\infty)$ 
 acts on $\frU(\infty)$ by left and right multiplication, thus we get a measure preserving action
of $\U(\infty)\times\U(\infty)$ on $\frU(\infty)$. 

More generally, we  fix $\lambda\in \C$ and  consider the probability measure
$$
\chi^n_\lambda = \prod_{k=1}^n\frac{\Gamma(k+\lambda)\Gamma(k+\ov\lambda)}
{\Gamma(k)\Gamma(k+\lambda+\ov\lambda)}
 \det(1+g)^\lambda \det(1+\ov g)^{\ov \lambda} d\chi^n(g)
$$
on $\U(n)$.
This system of measures also is projective and we get a family of measures $\chi_\lambda^\infty$
on the inverse limit of the chain (\ref{eq:chain}).
 
Applying the Cayley transform to matrices $g\in\U(n)$, we get the space of $n\times n$
Hermitian
matrices and the measures of the form
$$
C(\lambda, n)\det (1+iX)^{-\lambda-n}\cdot \det (1-iX)^{-\ov\lambda-n}
\,dX
.$$
Total integrals for such  measures (for real $\lambda$) were firstly evaluated by Hua Loo Keng
\cite{Hua}.

Now there is a substantial harmonic analysis on $\frU(\infty)$, see
\cite{Olsh-add2}, \cite{BO}.
Similar inverse limits  exists for all 10 series of compact semisimple symmetric spaces,
$$
\dots\longleftarrow G(n)/K(n)\longleftarrow  G(n+1)/K(n+1)\longleftarrow\dots
,$$
the group $G(\infty)$ acts on the inverse limit, see \cite{Pick}, \cite{Ner-hua}.
In \cite{Ner-hua} the group spaces $(K\times K)/K$, where $K=\SO(n)$, $\U(n)$, $\Sp(n)$
were considered. The Livshits map (\ref{eq:upsilon}) sends
  symmetric matrices ($g=g^t$) to symmetric,
  the space $\U(n)/\OO(n)$ of unitary symmetric matrices can be identified with 
 the real Lagrangian Grassmannian. Therefore we get a chain of Lagrangian Grassmannians
$$
\dots\longleftarrow \U(n)/\OO(n)\longleftarrow  \U(n+1)/\OO(n+1)\longleftarrow\dots
,$$
see \cite{Ner-gauss}, Section 3.6.
Other classical symmetric spaces can be reduced to the group case by the same trick
(we choose an appropriate involution in a group manifold commuting with the map $\Upsilon$).

There is also a similar construction for symmetric group (see \cite{KOV}), 
the limit object admits a substantial harmonic analysis. Our purpose is to obtain a $p$-adic counterpart of these
constructions.

\sm


{\bf\punct Notation.} Let

--- $\Q_p$ be the $p$-adic field;

--- $\bbO_p$ the ring of $p$-adic integers;

--- $\vol(\cdot)$ be the translation invariant $\sigma$-finite measure on
a linear space $\Q_p^m$ normalized by the condition
$\vol(\bbO_p^m)=1$;

--- $|\cdot|$ be the norm on $\Q_p$;

--- $\GL(n,\Q_p)$ and $\GL(n,\bbO_p)$ be the groups of invertible $n\times n$ matrices
over $\Q_p$ and $\bbO_p$;

--- $\Mat(n,\Q_p)$,  $\Mat(n,\bbO_p)$ be the spaces of all $n\times n$ matrices over 
$\Q_p$ and $\bbO_p$;

--- $\Symm(n,\Q_p)$, $\ASymm(n,\Q_p)$ be spaces of symmetric (skew-symmetric matrices) over $\Q_p$;

--- $g^t$ be transposed matrix;

--- $\Gr_{2n}^n$ is the Grassmannian of $n$-dimensional subspaces in $\Q_p^{2n}$.

\sm


{\bf\punct Measures $\mu_s^n$.}
Any $z\in\Mat(n,\Q_p)$ can be represented in the form
$$
z=A\begin{pmatrix}p^{-k_1}&0&\dots\\
0&p^{-k_2}&\dots\\
\vdots&\vdots&\ddots
 \end{pmatrix} B,
 $$
 where $A$, $B\in\GL(n,\bbO_p)$ and
 $$k_1\ge k_2\ge\dots\ge k_n\ge -\infty.
 $$
We say that $p^{k_j}$ are {\it singular numbers} of the matrix $z$.

We define the function $\gamma(z)$ on $\Mat(n,\bbO_p)$ by
$$
\gamma(z)=\prod_{k_j>0} p^{k_j}
$$ 
(we assume that a product of empty set of factors equals 1).

\begin{theorem}
\label{th:total-measure}
\begin{equation}
\int_{\Mat(n,\Q_p)} \gamma(z)^{-\alpha} \,d\vol(z)
=\prod_{j=1}^n\frac{1-p^{-\alpha+n-j}}{1-p^{-\alpha+n+j-1}}=:c(n,\alpha)
.
\label{eq:hua}
\end{equation}
The integral converges if $\alpha>2n-1$.
\end{theorem}

We define a measure $d\mu_s^n$, where $s>-1$, on $\Mat(n,\Q_p)$ by
$$
d\mu_s^n(z):=c(n,s+2n)^{-1}\gamma(z)^{-s-2n} d\vol_n(z)
.$$

The group $\GL(2n,\bbO_p)$ acts
on $\Mat(n,\Q_p)$ by linear-fractional transformations
\begin{equation}
z\mapsto (a+zc)^{-1}(b+zd)
,
\label{eq:linfrac}
\end{equation}
where $\begin{pmatrix}a&b\\c&d\end{pmatrix}\in \GL(2n,\bbO_p)$
is a block $(n+n)\times (n+n)$ matrix. Notice, that this formula
corresponds to the action of $\GL(2n,\bbO_p)$ on the Grassmannian
$\Gr_{2n}^n$.
Indeed, for an operator $z:\Q_p^n\to\Q_p^n$ consider its graph 
in $\Q_p^n\oplus\Q_p^n$, it consists of vectors
$$
v\oplus vz, \,\,\text{where $v\in\Q_p^n$ is a row matrix.}
$$
The get a chart on Grassmannian, the complement of the chart has zero measure.
A verification of (\ref{eq:linfrac}) is straightforward (see, e.g., \cite{Ner-gauss}, Theorem 2.3.1).

\begin{theorem}
\label{th:derivative}
{\rm a)} For any $\begin{pmatrix}a&b\\c&d\end{pmatrix}\in\GL(2n,\Q_p)$,
\begin{equation}
d\mu_s^n\left((a+zc)^{-1}(b+zd)\right)=|\det(a+zc)|^{s} d\mu_s^n(z).
\label{eq:derivative}
\end{equation}

{\rm b)} The measure $\mu_s^n$ is a unique probability Borel measure
on $\Mat(n,\Q_p)$ satisfying the equation
\begin{equation}
d\nu\left((a+zc)^{-1}(b+zd)\right)=|\det(a+zc)|^{s} d\nu(z)
\label{eq:derivative-2}
\end{equation}
 for any $\begin{pmatrix}a&b\\c&d\end{pmatrix}\in \GL(2n,\bbO_p)$.
\end{theorem}

In particular, for $s=0$ we get a unique $\GL(2n,\Q_p)$-invariant measure on the
Grassmannian $\Gr_{2n}^n$.

\sm


{\bf\punct Projective limits.} 
Consider a $(n+1)\times(n+1)$ matrix
$$z=\begin{pmatrix}
z_{11}&z_{12}\\
z_{21}&z_{22}
\end{pmatrix}
.$$
 Consider the map $\Pi:\Mat(n+1,\Q_p)\to \Mat(n,\Q_p)$
given by
\begin{equation}
\Pi: \begin{pmatrix}
z_{11}&z_{12}\\
z_{21}&z_{22}
\end{pmatrix}\mapsto z_{11}.
\label{eq:Pi}
\end{equation}

\begin{theorem}
\label{th:projective}
The pushforward of the measure $\mu_s^{n+1}$ under the map $\Pi$ is $\mu_s^{n}$
\end{theorem}

Thus we get a chain 
$$\dots
\longleftarrow \bigl(\Mat(n,\Q_p),d\mu_s^n\bigr)
\longleftarrow \bigl(\Mat(n+1,\Q_p),d\mu_s^{n+1}\bigr)\longleftarrow\dots
$$
By the Kolmogorov theorem (see, e.g. \cite{Shi}, Section 2.3, Theorem 3)
 the inverse limit in the category of measure spaces
 is well defined, denote by $d\mu_s^\infty$   the inverse limit of measures $\mu_s^n$,
  this measure can be regarded
 as a measure on $\Mat(\infty,\Q_p)\simeq \Q_p^{\infty\times\infty}$.
 
 
 \sm
 
 {\bf\punct Symmetries of the measures $\mu_s^\infty$.}
 Consider the chain of groups
 $$
 \dots \longrightarrow \GL(2n,\bbO_p) \longrightarrow
 \GL\bigl(2(n+1),\bbO_p\bigr)  \longrightarrow \dots
 $$
 and its inductive limit $\GL(2\infty,\bbO_p)$. In other words, $\GL(2\infty,\bbO_p)$
 is the group of $(\infty+\infty)\times (\infty+\infty)$ matrices $g$ with integer elements
 such that $g^{-1}$ also has integer elements and $g-1$ has only finite number of non-zero entries.
 
 \begin{proposition}
 \label{pr:quasiinvariance}
{\rm a)} The measure $\mu_s$ is quasiinvariant with respect to the action
 $z\mapsto (a+zc)^{-1} (b+zd)$ of $\GL(2\infty,\bbO_p)$.
  The Radon--Nikodym derivative is
 $|\det(a+zc)|^{s}$.
 
 \sm
 
 {\rm b)} In particular for $s=0$ the measure $\mu_s^\infty$ is $\GL(2\infty,\Q_p)$-invariant.
 
 \sm
 
 {\rm c)} The measure $\mu_s^\infty$ is invariant with respect to the subgroup
 $P\subset \GL(2\infty,\Q_p)$ consisting
 of matrices $\begin{pmatrix}a&b\\0&d\end{pmatrix}$.
  \end{proposition}
 
 We must explain the meaning of the expression $|\det(a+zc)|^{s}$.
 Represent $g\in \GL(2\infty,\bbO_p)$ as
 a block matrix of size $(k+\infty+k+\infty)$, where $k$ is sufficiently large,
 $$
g= \begin{pmatrix}a&b\\c&d\end{pmatrix}=
\begin{pmatrix}a_{11}&0&b_{11}&0\\
0&1&0&0\\
c_{11}&0&d_{11}&0
\\
0&0&0&1\end{pmatrix}
.
 $$
 Represent $z\in\Mat(\infty,\Q_p)$ as a block matrix  of size $(k+\infty)$:
 $$
 z=\begin{pmatrix}
z_{11}&z_{12}\\
z_{21}&z_{22}
\end{pmatrix}
 .$$
 Then
 \begin{multline*}
|\det(a+zc)|=\left|\det\begin{pmatrix} 
 a_{11}&0\\0&1
 \end{pmatrix}
 +
 \begin{pmatrix}
z_{11}&z_{12}\\
z_{21}&z_{22}
\end{pmatrix}
\begin{pmatrix}
c_{11}&0\\0&0
\end{pmatrix}\right|
=\\=\left|
\det\begin{pmatrix} 
 a_{11}+z_{11}c_{11}&0\\z_{12}c_{11}&1
 \end{pmatrix}
 \right|=\bigl|\det  (a_{11}+z_{11}c_{11})\bigr|
, \end{multline*}
and we get a determinant of a finite  matrix.

\smallskip

Next, we define two completions 
$$
\ov{\ov{\GL}}(2\infty,\bbO_p)\supset {\ov{\GL}}(2\infty,\bbO_p)\supset {{\GL}}(2\infty,\bbO_p)
$$
of the group $\GL(2\infty, \bbO_p)$.

First, consider the group $\ov{\ov T}$   consisting of $(\infty+\infty)\times (\infty+\infty)$
matrices over $\bbO_p$ having the form
\begin{equation}
h=
\begin{pmatrix}a&b\\0&d\end{pmatrix},
\label{eq:T1}
\end{equation}
where
\begin{equation}
a=\begin{pmatrix} \alpha_{11}&0&0&\dots\\
\alpha_{21}&\alpha_{22}&0&\dots\\
\alpha_{31}&\alpha_{32}&\alpha_{33}&\dots\\
\vdots&\vdots&\vdots&\ddots
\end{pmatrix},
\qquad
d=
\begin{pmatrix}
\delta_{11}&\delta_{12}&\delta_{13}&\dots\\
0&\delta_{22}&\delta_{23}&\dots\\
0&0&\delta_{33}&\dots\\
\vdots&\vdots&\vdots&\ddots
\end{pmatrix},
\label{eq:T2}
\end{equation} 
and
\begin{equation}
|\alpha_{11}|=|\alpha_{22}|=\dots=1\qquad |\delta_{11}|=|\delta_{22}|=\dots=1
\label{eq:T3}
.
\end{equation}
As a set the group $\ov{\ov T}$ is a direct product of countable number of copies
of $\bbO_p\setminus p\bbO_p$  (corresponding to $\alpha_{jj}$ and $\delta_{kk}$)
and  countable number of copies of $\bbO_p$ (corresponding to remaining $\alpha_{ij}$,
$\delta_{kl}$, and to $\beta_{mn}$).

Also consider a smaller group $\ov T\subset\ov{\ov T}$ consisting of the matrices having the form
(\ref{eq:T1})--(\ref{eq:T2}) but (\ref{eq:T3}) is replaced by
\begin{equation}
\alpha_{11}=\alpha_{22}=\dots=1\qquad \delta_{11}=\delta_{22}=\dots=1
.
\label{eq:T4}
\end{equation}

{\sc Remark.}
Notice that the matrices $h\in \ov{ \ov T}$ can be made upper triangular after a permutation
of  basis elements (the first $\infty$ of basis elements must be written in the inverse
order). \hfill $\square$

\sm

Denote by $\ov{\ov\GL}(2\infty,\bbO_p)$ the group of matrices generated by
$\GL(2\infty,\bbO_p)$ and $\ov{\ov T}$. This group consists of invertible matrices
$g=\begin{pmatrix}a&b\\c&d\end{pmatrix}$ over $\bbO_p$
such that

$1^\circ$. $c$ has only finite number nonzero entries;

$2^\circ$. $a$ has only finite number of nonzero entries upper the diagonal;

$3^\circ$. $d$ has only finite number of nonzero entries lower the diagonal;

$4^\circ$. $g^{-1}$ has integer elements.

\sm

Denote by $\ov\GL(2\infty,\bbO_p)\subset\ov{\ov\GL}(2\infty,\bbO_p)$
the group of matrices generated by
$\GL(2\infty,\bbO_p)$ and $\ov T$. We must replace $2^\circ$ and $3^\circ$ by:

$2^{\circ\circ}$. $a-1$ has only finite number of nonzero entries on the diagonal and upper the diagonal;

$3^{\circ\circ}$. $d-1$ has only finite number of nonzero entries on the diagonal 
and lower the diagonal.

\begin{theorem}
\label{th:quasiinvariance}
The group $\ov{\ov\GL}(2\infty,\bbO_p)$ acts on $\Mat(\infty,\Q_p)$
by transformations
 $z\mapsto (a+zc)^{-1} (b+zd)$ leaving the measure $\mu_s^\infty$ quasiinvariant.
  The Radon--Nikodym derivative is
 $|\det(a+zc)|^{s}$.
\end{theorem}

The meaning of the expressions $(a+zc)^{-1} (b+zd)$ and 
$|\det(a+zc)|^{s}$ 
will be explained in Subsection \ref{ss:big-group}.

\begin{proposition}
\label{pr:for-induction}
For any $g\in {\ov\GL}(2\infty,\bbO_p)$ and $z\in \Mat(\infty,\Q_p)$, the expression
$$
\det (a+zc)\in \Q_p
$$
is well defined.
\end{proposition}

Denote by $\Q_p^\times$, $\C^\times$ the multiplicative groups of $\Q_p$ and $\C$.
Let $\chi$ be a homomorphism $\Q_p^\times\to\C^\times$, $|\chi(z)|=1$.
Then we can define a unitary representation of ${\ov\GL}(2\infty,\bbO_p)$
in $L^2\bigl(\Mat(\infty,\Q_p),\mu_s^\infty\bigr)$ by the formula
\begin{equation}
\rho_{s,\chi} \begin{pmatrix}a&b\\c&d\end{pmatrix} f(z) = f\bigl((a+zc)^{-1} (b+zd)\bigr)
\det(a+zc)^{s/2}\chi(a+zc)
.
\label{eq:representation}
\end{equation}

\sm


{\bf \punct Some remarks on real-$p$-adic parallel.}
Analogs of noncompact Riemannian symmetric spaces over $p$-adic numbers are Bruhat--Tits buildings
(here there is a well-known and deep parallel, see a discussion and further references in
 \cite{Ner-gauss}). Analogs of Hua
integrals for buildings exist (see \cite{Ner-bui}) and  they are used below in Section 2. But it seems
that they do not admit projective limits (moreover, real-$p$-adic  analogy does not requires this,
for noncompact Riemannian symmetric spaces there are no
inverse limits, see 
\cite{Ors}). 

Apparently,
 there are no reasonable $p$-adic analogs of compact Riemannian symmetric spaces. 
 
 On the other hand,
 the classical compact Riemannian symmetric spaces  are Grassmannians or isotropic Grassmannians 
 (see, e.g., \cite{Ner-gauss}, Subsections D.1), our construction is an emulation of Pickrell's approach \cite{Pick}.


\sm

{\bf\punct   Other inverse limits of Grassmannians.} 

a) {\it Symplectic Lagrangian Grassmannian.} Consider the group $\Sp(2n,\bbO_p)$ consisting of
integer
$(n+n)\times(n+n)$ matrices $\begin{pmatrix} a&b\\c&d\end{pmatrix}$
preserving skew-symmetric bilinear form $\begin{pmatrix}0&1\\-1&0 \end{pmatrix}$  in $\Q_p^n\oplus \Q_p^n$.
Consider the Grassmannian $\mathrm{L}_n$ of Lagrangian subspaces in $\Q_p^{2n}$.
Almost all elements of $\mathrm{L}_n$  can be represented as graphs of operators
$\Q_p^n \oplus 0\to 0\oplus \Q_p^n$, the corresponding matrices $z$ are symmetric (see, e.g., 
\cite{Ner-gauss}, Theorem 3.1.4). The action of
the group $\Sp(2n,\bbO_p)$ on $\Symm(n,\Q_p)$ is given by the same formula (\ref{eq:linfrac}).

We define measures
\begin{equation}
d\mu_s^n(z)=\gamma(z)^{-s-n-1}d\vol_n(z)
\label{eq:mu-symplectic}
\end{equation}
on $\Symm(n,\Q_p)$, we choose normalizing constants  $c(s,n)$
to obtain probability measures.
 Their Radon--Nikodym derivatives with respect to transformations
 in $\Sp(2n,\bbO_p)$ are
\begin{equation}
\frac{d\mu_s^n\bigl((a+zc)^{-1}(b+zd) \bigr)}{d\mu_s^n(z)}
=|\det(a+zc)|^s
.
\label{eq:RN}
\end{equation}

Next, these measures
 form a projective system 
with respect to the map $\Pi:\Symm(n+1,\Q_p)\to \Symm(n,\Q_p)$, see (\ref{eq:Pi}).

Now we can consider the inverse limit of measure spaces
$$
\dots \longleftarrow \bigl(\Symm(n,\Q_p),d\mu_s^n\bigr)
\longleftarrow \bigl(\Symm(n+1,\Q_p),d\mu_s^{n+1}\bigr)\longleftarrow \dots
$$

b) {\it Isotropic orthogonal Grassmannians.} Consider the group $\OO(2n,\bbO_p)$ consisting of
integer
$(n+n)\times(n+n)$ matrices $\begin{pmatrix} a&b\\c&d\end{pmatrix}$
preserving the symmetric bilinear form $\begin{pmatrix}0&1\\1&0 \end{pmatrix}$
in $\Q_p^n\oplus \Q_p^n$. Consider the subgroup $\SO(2n,\bbO_p)\subset\OO(2n,\bbO_p)$
consisting of matrices with determinant $=1$.
Consider the Grassmasnnian $\mathrm{Is}_n$ of $n$-dimensional isotropic subspaces.
The group $\SO(2n,\bbO_p)$ has two orbits on $\mathrm{Is}_n$, we choose
one. Namely,
 consider
the set $\mathrm{Is}_n^0$ of $n$-dimensional isotropic
subspaces $M$ such that $M\cap (0\oplus\Q_p^n)$ has even dimension,
see, e.g.,
\cite{Ner-book}, Proposition 2.2.2, Lemma 3.3.1. Almost all elements of 
$\mathrm{Is}_n^0$ are graph of operators $z:\Q_p^n\oplus 0\to 0\oplus \Q_p^n$
and matrices $z$ are skew-symmetric ($z=-z^t$). 
We define probability measures
\begin{equation}
d\mu_s^n(z)= a(s,n) \gamma(z)^{-s-n+1}d\vol_n(z).
\label{eq:mu-orthogonal}
\end{equation}
on $\ASymm(n,\Q_p)$. Their Radon--Nikodym derivatives
are given by the same formula (\ref{eq:RN}).
Measures $d\mu_s^n$ form a projective system with respect to the maps
 $\Pi:\ASymm(n+1,\Q_p)\to \ASymm(n,\Q_p)$.
Again, we can consider  inverse limits of measure spaces
$$
\dots \longleftarrow \bigl(\ASymm(n,\Q_p),d\mu_s^n\bigr)\longleftarrow 
\bigl(\ASymm(n+1,\Q_p,d\mu_s^{n+1}\bigr),\longleftarrow \dots
$$

c) In both cases (symplectic and orthogonal), the situation is parallel to the picture described above.

But the author does not know the explicit analog of the formula (\ref{eq:hua})
for complete measure, because our prove is 
a reduction to the beta-function of Bruhat--Tits buildings \cite{Ner-bui}, which was evaluated only
for  $\GL(n,\Q_p)$-case. Also, we can normalize the measures (\ref{eq:mu-symplectic}), 
(\ref{eq:mu-orthogonal})
only if they are finite. Evidently, this holds for $s=0$ since in this case
we have measure on Grassmannians invariant with respect to transitive actions
 compact groups ($\Sp(2n,\bbO_p)$ and $\SO(2n, \bbO_p)$). Since $\gamma(z)^{-\alpha}$
 decreases as a function $\alpha$, we get that our measures are well-defined
 at least for $s\ge 0$. 

\sm


{\bf\punct Remark. On representations of infinite-dimensional classical  $p$-adic groups.}
Basic representation theory of infinite-dimensional classical groups 
and infinite symmetric groups was developed in
70-80s, see \cite{VK1}, \cite{VK2}, \cite{Olsh-symm} for symmetric groups
and \cite{Olsh-short}, \cite{OlshGB}, \cite{Olsh-topics},
 \cite{Pick2}, \cite{Olsh-topics}, \cite{Ner-book} for classical groups. 
These works had various continuations, see, e.g., \cite{Ner-book}, \cite{BO}, \cite{Olsh-add2},
\cite{KOV}, and further  references in \cite{Ner-char}.

Representations of infinite-dimensional classical $p$-adic groups remains to be a non well-understood topic. Now   two substantial constructions are known.
The first is the Weil representation of infinite-dimensional symplectic group
and the corresponding contractive semigroup
 (Nazarov%
 \footnote{A weaker version of construction  is in \cite{Zel}.},
  \cite{NNO}, \cite{Naz}, and  a partial exposition in  \cite{Ner-gauss}, Sections 10.7, 11.2). 
 The second is the multiplication of double cosets and
a $p$-adic analog of 
characteristic operator-function, see \cite{Ner-char}. 

For infinite-dimensional groups over  finite fields, see
 \cite{VK3}. 


\sm

{\bf\punct Further structure of the paper.}
In Section 2 we evaluate the integral (\ref{eq:hua}).
In Section 3 we prove Theorem \ref{eq:derivative} on transformations
of measures $\mu_s^n$. In Section 4 we prove the statements about 
measures $\mu_s^\infty$.



\section{Proofs. Calculation of the integral}

\COUNTERS

Here we prove Theorem \ref{th:total-measure}.
In this section,
$$
G:=\GL(n,\Q_p),\qquad K:=\GL(n,\bbO_p).
$$
Denote by $\Lat(n)$ the set of lattices  (see, e.g., \cite{Ner-gauss}, Section 10.3)
 in $\Q_p^n$, we have
$$
\Lat(n)\simeq G/K.
$$

\sm

{\bf\punct Properties of the function $\gamma$.} The following statement
is obvious.

\begin{lemma}
 Let $z\in\Mat(n,\Q_p)$, $\det(z)\ne 0$. Then
\begin{align}
\frac{\gamma(z)}{\gamma(z^{-1})}
&=|\det(z)|;
\\
\gamma(z)&=\vol(z\bbO_p^n+\bbO_p^n);
\\
\gamma(z^{-1})^{-1}&=\vol(z\bbO_p^n\cap\bbO_p^n).
\end{align}
\end{lemma}

Also, note that 
\begin{equation}
|\det(z)|=\vol(z\bbO_p^n).
\end{equation}


{\bf\punct  Haar measure.} According \cite{Mac}, Section V.2,
the Haar measure on $G$ is given by
\begin{equation}
|\det(z)|^{-n}\,d\vol(z).
\label{eq:haar}
\end{equation}
It is convenient to normalize this measure as
$$
d\chi(z)=\frac 1{\vol (K)}
|\det(z)|^{-n}\,d\vol(z),
$$
then $\chi(K)=1$.

\begin{lemma}
$$
\vol(K)= \prod_{j=1}^n\left(1-p^{-j}\right).
$$
\end{lemma}

{\sc Proof.} Consider the natural map $\Mat(n,\bbO_p)\to \Mat(n,\bbF_p)$, 
where $\bbF_p=\bbO_p/p\bbO_p$ is the field with $p$ elements. The total number
of points in $\Mat(n,\bbF_p)$ is $p^{n^2}$. The total number of points
in  $\GL(n, \bbF_p)\subset \Mat(n,\bbF_p)$ is  $\prod_{j=1}^n (p^n-p^j)$.
\hfill $\square$

\sm

{\bf\punct  Calculation of the integral.}
Keeping in the mind the expression for Haar measure we transform our integral
(\ref{eq:hua}) as
\begin{multline*}
\int\limits_{\Mat(n,\Q_p)} \gamma(z)^{-t}\, d\vol(z)=
\vol(K)
\int\limits_{G} |\det(z)|^{n}\, \gamma(z)^{-t}\,  d\chi(z)
=\\=
\vol(K)
\int\limits_{G} |\det(z)|^{-t+n}
 \vol(z\bbO_p^n\cap \bbO_p^n)^{t}   d\chi(z)
 =\\=
\vol(K)
\int\limits_{G} \vol(z\bbO_p^n)^{-t+n} \,
 \vol(z\bbO_p^n\cap \bbO_p^n)^{t} \,  d\chi(z).
\end{multline*}
The integrand is constant on each coset $zK\subset G$, by the invariance
of Haar measure, we have $\chi(zK)=1$. Therefore we come to a summation over
$G/K=\Lat(n)$:
$$
\vol(K)
\sum_{Q\in\Lat(n)}
\vol(Q)^{-t+n}
 \vol(Q\cap \bbO_p^n)^{t} 
 .
$$
This expression is a special case of the 'beta-function of Bruhat--Tits building' evaluated
in \cite{Ner-bui}, Theorem 2.1. In notation of \cite{Ner-bui}, we set
$$\alpha_1=\dots=\alpha_n=t\qquad \beta_1=\dots=\beta_n=-t+n$$
and get
$$
\vol(K)\prod_{j=1}^n
\frac{1-p^{-t+n-j}}{(1-p^{-t+n+j-1})(1-p^{-j})}=
\prod_{j=1}^n
\frac{1-p^{-t+n-j}}{(1-p^{-t+n+j-1})}
.$$


\section{Proofs. Transformation of measures}

Here we prove Theorem \ref{th:quasiinvariance}.
 
\COUNTERS

{\bf\punct Formula for Radon--Nikodym derivative.} 1) Consider the subgroup 
$P\subset\GL(2n,\bbO_p)$
 consisting of matrices  $\begin{pmatrix}a&b\\0&d \end{pmatrix}$,
it acts on $\Mat(n,\Q_p)$ by the transformations
$$
z\mapsto a^{-1} (b+zd)= a^{-1} b+a^{-1}zd.
$$
The measure $d\vol(z)$ is invariant with respect to such transformations. The function
$\gamma(z)^{-s-2n}$ also is invariant. Therefore the measure $\mu_s^n$ is $P$-invariant.

\sm

2) Consider the matrix $\begin{pmatrix}0&1\\1&0\end{pmatrix}$.
The corresponding transformation is $z\mapsto z^{-1}$.
We must show that
\begin{equation}
\gamma(z^{-1})^{-s-2n} d\vol(z^{-1})= |\det(z)|^{s}\gamma(z)^{-s-2n}d\vol(z) 
\label{eq:1}
.
\end{equation}
First 
$$d\vol(z^{-1})=|\det(z)|^{-2n}d\vol(z) ,$$
this is equivalent to the formula (\ref{eq:haar}) for the Haar measure.
On the other hand,
$$
\gamma(z^{-1})^{-s-2n}=\gamma(z)^{-s-2n}\,|\det z|^{s+2n}
,$$
and we get (\ref{eq:1}).

\sm

3) It can be easily shown that $\GL(2n,\bbO_p)$ is generated by subgroup
$P$ and $\begin{pmatrix}0&1\\1&0 \end{pmatrix}$.
Denote 
\begin{equation}
z*g:=(a+zc)^{-1}(b+zd)
.\end{equation}
The Radon--Nikodym derivative 
$$
c(g,z):=\frac{d\mu_s^n(z*g)}{d\mu_s^n(z)}
.
$$
satisfies the chain rule
\begin{equation}
c(g_1g_2,z)=c(g_1,z)\,c(g_2,z*g_1)
.
\label{eq:chain-rule}
\end{equation}
On the other hand, the expression 
$$
\wt c(g,z):=|\det(a+zc)|^{s}
$$
also satisfies the chain rule. Since $c(h,z)=\wt c(h,z)$ for generators of $G$,
they coincide everywhere.

\sm


{\bf\punct Uniqueness of the measure.}
Consider a measure $\nu$ satisfying  equation (\ref{eq:derivative-2}).
First, we consider matrices $\begin{pmatrix}1&b\\0&1\end{pmatrix}$
and get that $\nu$ is invariant with respect to translations
$$
z\mapsto z+b,\qquad b\in\Mat(n,\bbO_p)
.$$
Therefore $\nu$ has the form
$$
d\nu(z)=f(z)\,d\mu_s^n(z)
,$$
where $f(z)$ is a locally constant function.
Evidently, $f(z)$ is $\GL(2n,\bbO_p)$-invariant. Since the action of $\GL(2n,\bbO_p)$ on the Grassmannian
is transitive, we get that $f(z)$ is constant.


\section{Proofs. Projective limits}

\COUNTERS

{\bf\punct Proof of Theorem \ref{th:projective}.}
It is sufficient to show that the $\Pi$-pushforward $\nu$ of $\mu_s^{n+1}$ satisfies
the quasiinvariance property (\ref{eq:derivative-2}).
It is sufficient to verify this property for generators 
of $\GL(2n,\bbO_p^n)$:
$$
g_1=
\begin{pmatrix}a&0\\ 0&d\end{pmatrix}, \qquad g_2=
\begin{pmatrix}1&b\\ 0&1\end{pmatrix},\qquad g_3=
\begin{pmatrix}0&1\\ 1&0\end{pmatrix}
.$$

We consider the corresponding elements of $\GL(2(n+1),\bbO_p)$:
$$
\wt g_1=
\begin{pmatrix}a&0&0&0\\
0&1&0&0\\
0& 0&d&0\\
0&0&0&1
\end{pmatrix}
, \qquad \wt g_2=
\begin{pmatrix}1&0&b&0\\
0&1&0&0\\
 0&0&1&0\\
0&0&0&1 
 \end{pmatrix},\qquad \wt g_3=
\begin{pmatrix}0&0&1&0\\0&1&0&0\\ 1&0&0&0\\
0&0&0&1
\end{pmatrix}
$$
(sizes of the matrices are $n+1+n+1$). We have
\begin{align}
\wt g_1:
 \begin{pmatrix}
z_{11}&z_{12}\\
z_{21}&z_{22}
\end{pmatrix} &\mapsto
\begin{pmatrix}
a&0\\0&1
\end{pmatrix}^{-1}
 \begin{pmatrix}
z_{11}&z_{12}\\
z_{21}&z_{22}
\end{pmatrix}
\begin{pmatrix}
d&0\\0&1
\end{pmatrix}=
\begin{pmatrix}
a^{-1} z_{11}d & *\\ *&*
\end{pmatrix}
;\\
\wt g_2:
 \begin{pmatrix}
z_{11}&z_{12}\\
z_{21}&z_{22}
\end{pmatrix}&\mapsto
\begin{pmatrix}
b&0\\0&1
\end{pmatrix}+
 \begin{pmatrix}
z_{11}&z_{12}\\
z_{21}&z_{22}
\end{pmatrix}= 
\begin{pmatrix}
b+z_{11}&*\\
*&*
\end{pmatrix}
;\end{align}
\begin{multline*}
\wt g_3:
 \begin{pmatrix}
z_{11}&z_{12}\\
z_{21}&z_{22}
\end{pmatrix}\mapsto
\\\mapsto
\left[
\begin{pmatrix}
0&0\\0&1
\end{pmatrix}+
\begin{pmatrix}
1&0\\0&0
\end{pmatrix}
 \begin{pmatrix}
z_{11}&z_{12}\\
z_{21}&z_{22}
\end{pmatrix}\right]^{-1}
\left[
\begin{pmatrix}
1&0\\0&0
\end{pmatrix}+
\begin{pmatrix}
0&0\\0&1
\end{pmatrix}
 \begin{pmatrix}
z_{11}&z_{12}\\
z_{21}&z_{22}
\end{pmatrix}\right]
=\\=
 \begin{pmatrix}
z_{11}&0\\
z_{21}&1
\end{pmatrix}^{-1}
 \begin{pmatrix}
1&0\\
z_{21}&z_{22}
\end{pmatrix}
= \begin{pmatrix}
z_{11}^{-1}&*\\
*&*
\end{pmatrix}
.\end{multline*}

 Thus we get
$$
\Pi(z*\wt g_j)=(\Pi z)*g_j.
$$
In the first two cases the Radon--Nikodym derivative is 1, in the last case
$|\det z_{11}|^{s}$.
Thus,
$$
\frac{d\mu_s^{n+1}(z*\wt g_j)}{d\mu_s^{n+1}(z)}=\frac{d\mu_s^{n}((\Pi z)* g_j)}{d\mu_s^{n}(\Pi z)}
.$$
This implies desired property of the $\Pi$-pushforward of $\mu_s^{n+1}$.


\sm

{\bf\punct An abstract lemma.}
Let $(\Omega_j,\mu_j)$ be a sequence of  Lebesgue measure  spaces with
  probability measures. Let $\Pi_k^j:\Omega_j\to \Omega_k$, where $k<j$, be maps
  such that $\Pi^k_l\Pi_k^j=\Pi^j_l$ and the $\Pi^j_k$-pushforward of 
  $\mu_j$ is  $\mu_k$. Denote by $(\Omega_\infty,\mu_\infty)$ the
  projective limit of the chain
  $$
  \dots\longleftarrow \Omega_k \longleftarrow  \Omega_{k+1} \longleftarrow \dots
  $$

  Let $G$ be a group. Let for any $g\in G$ there exists
  $j$ such that for all $k\ge j$ there is a  transformation
  $g_{[k]}:\Omega_k\to\Omega_k$ leaving the measure $\mu_j$ quasiinvariant and for
  $l>k>j$ we have (a.s.)
  \begin{equation}
  \Pi_k^l \bigl( g_{[l]}( \omega)\bigr) =
  g_{[k]} \bigl( \Pi_k^l (\omega)\bigr), \qquad \omega\in\Omega_l,
  \label{eq:comp1}
  \end{equation}
  and (a.s.)
  \begin{equation}
  (g_{[k]})'(\Pi_{[k]}^l\omega)=(g_{[l]})'(\omega), \qquad \omega\in\Omega_l.
   \label{eq:comp2}
  \end{equation}
  Let also for any $g$, $h\in G$ for sufficiently large $m$
  $$
  (gh)_{[m]}= g_{[m]} h_{[m]}
  .
  $$
 
 \begin{lemma} 
 \label{l:l}
  Under these conditions there is action of $G$ on $\Omega_\infty$ by transformations
  $g_{[\infty]}$
  leaving the measure $\mu_\infty$ quasiinvariant, they are determined by
  $$
  \Pi^\infty_k \bigl( g_{[\infty]} (\omega)\bigr)=
    g_{[k]}\bigl( \Pi^\infty_k( \omega\bigr), \qquad \omega\in\Omega_\infty,
  $$
  and
  $$
  (g_{[k]})'(\Pi_{[k]}^\infty\omega)=(g_{[\infty]})'(\omega), \qquad \omega\in\Omega_\infty.
  $$
      \end{lemma}

This is straightforward.

 Proposition \ref{pr:quasiinvariance} is an immediate corollary of the lemma.

\sm


{\bf \punct Action of $\ov{\ov \GL}(2\infty,\bbO_p)$.%
\label{ss:big-group}} We wish to reduce
 Theorem \ref{th:quasiinvariance} to Lemma \ref{l:l}.
 Now $\Omega_k:=\Mat(k,\Q_p)$, projections $\Omega_l\to\Omega_k$ are cutting
 of left upper $k\times k$ corner. We must construct transformations
  $g_{[k]}:\Mat(k,\Q_p)\to \Mat(k,\Q_p)$.
 
 Fix $g\in \ov\GL(2\infty,\bbO_p)$.
Choose a sufficiently large $k$ and represent 
$g=\begin{pmatrix}a&b\\c&d\end{pmatrix}$
in the form
$$
g=
\begin{pmatrix}a_{11}&0&b_{11}&b_{12}\\
a_{21}&a_{22}&b_{21}&b_{22}\\
c_{11}&0&d_{11}&d_{12}
\\
0&0&0&d_{22}\end{pmatrix}
=
\begin{pmatrix}a_{11}^{(k)}&0&b_{11}^{(k)}&b_{12}^{(k)}\\
a_{21}^{(k)}&a_{22}^{(k)}&b_{21}^{(k)}&b_{22}^{(k)}\\
c_{11}^{(k)}&0&d_{11}^{(k)}&d_{12}^{(k)}
\\
0&0&0&d_{22}^{(k)}\end{pmatrix}
$$
(below we sometimes omit upper index $^{(k)}$).
We formally calculate 
\begin{multline*}
(a+zc)^{-1} (b+zd)=
\\
\left[
\begin{pmatrix}a_{11}&0\\a_{21}&a_{22}\end{pmatrix}
+
\begin{pmatrix}z_{11}&z_{12}\\z_{21}&z_{22}\end{pmatrix}
\begin{pmatrix}
c_{11}&0\\0&0
\end{pmatrix}
\right]^{-1}
\left[
\begin{pmatrix}b_{11}&b_{12}\\b_{21}&b_{22}\end{pmatrix}
+\begin{pmatrix}z_{11}&z_{12}\\z_{21}&z_{22}\end{pmatrix}
\begin{pmatrix} d_{11}&d_{12}\\ 0&d_{22} \end{pmatrix}
\right]
\\
=
\begin{pmatrix} (a_{11}+z_{11}c_{11})^{-1}&0\\ *&* \end{pmatrix}
\begin{pmatrix} b_{11}+z_{11}d_{11}&*\\ *&* \end{pmatrix}
=\\=
\begin{pmatrix} (a_{11}^{(k)}+z_{11}^{(k)}c_{11}^{(k)})^{-1}
(b_{11}^{(k)}+z_{11}^{(k)}d_{11}^{(k)})&*\\ *&* \end{pmatrix}
.
\end{multline*}
We observe that for sufficiently large $k$
the $k\times k$ left upper corner of $(a+zc)^{-1} (b+zd)$
depends only on the $z_{11}^{(k)}$ . Now we assign the transformation
$$
g_{[k]}: u
\mapsto
(a_{11}^{(k)}+uc_{11}^{(k)})^{-1}(b_{11}^{(k)}+ud_{11}^{(k)})
$$ 
of $\Mat(k,\Q_p)$.
The same calculation shows the compatibility (\ref{eq:comp1}).

Next, let us write formally the Radon--Nikodym derivative of our transformation 
is 
\begin{equation}
\left| \det
\begin{pmatrix} a_{11}+z_{11}c_{11}&0\\ a_{21}+z_{21}c_{11} &a_{22} \end{pmatrix}
\right|^{s}=|\det (a_{11}+z_{11}c_{11})|^{s}|\det a_{22}|^{s}
\label{eq:last}
.
\end{equation}
But $a_{22}$ is lower triangular and we can set $|\det a_{22}|=1$.

Note that the expression(\ref{eq:last})  coincides with the Radon--Nikodym derivative of $g_{[k]}$
and does not change under a pass $k\to k+1$. Therefore we have compatibility
(\ref{eq:comp1}).


\sm

{\bf \punct Proof of Proposition \ref{pr:for-induction}.}
For fixed $g=\begin{pmatrix}a&b\\c&d\end{pmatrix}\in\ov\GL(2n,\Q_p)$,
The following expression 
$$f(z)=
 \det
\begin{pmatrix} a_{11}+z_{11}c_{11}&0\\ a_{21}+z_{21}c_{11} &a_{22} \end{pmatrix}
:=\det (a_{11}+z_{11}c_{11})
$$
is
a well-defined function on $\Mat(\infty,\Q_p)$.

Also note that this function satisfies the chain rule
(\ref{eq:chain-rule}). Therefore formula (\ref{eq:representation})
determines a representation of the group $\ov{\GL}(2\infty,\bbO)$.

\tt Math.Dept., University of Vienna,

 Nordbergstrasse, 15,
Vienna, Austria

\&

Institute for Theoretical and Experimental Physics,

Bolshaya Cheremushkinskaya, 25, Moscow 117259,
Russia

\&

Mech.Math.Dept., Moscow State University,

Vorob'evy Gory, Moscow

e-mail: neretin(at) mccme.ru

URL:www.mat.univie.ac.at/$\sim$neretin

wwwth.itep.ru/$\sim$neretin

\end{document}